\documentclass{amsart}

\usepackage{amsmath,amsfonts,amssymb,amscd,verbatim,delarray,url}

\font\cyr=wncyr8 scaled \magstep1%

\def\Sha{\text{\cyr Sh}}

\newcommand{\F}{{\mathbb F}}

\newcommand{\Q}{{\mathbb Q}}

\DeclareMathOperator{\Inn}{Inn}
\DeclareMathOperator{\End}{End}
\DeclareMathOperator{\Aut}{Aut}
\DeclareMathOperator{\Endc}{End_{\text{\rm{c}}}}
\DeclareMathOperator{\Autc}{Aut_{\text{\rm{c}}}}
\DeclareMathOperator{\Out}{Out}
\DeclareMathOperator{\Outc}{Out_{\text{\rm{c}}}}

\DeclareMathOperator{\FSym}{FSym} \DeclareMathOperator{\FAlt}{FAlt}
\DeclareMathOperator{\Sym}{Sym} \DeclareMathOperator{\Alt}{Alt}

\begin{document}
\numberwithin{equation}{section}

\newtheorem{theorem}{Theorem}[section]
\newtheorem{lemma}[theorem]{Lemma}

\newtheorem{prop}[theorem]{Proposition}
\newtheorem{proposition}[theorem]{Proposition}
\newtheorem{corollary}[theorem]{Corollary}
\newtheorem{corol}[theorem]{Corollary}
\newtheorem{conj}[theorem]{Conjecture}
\newtheorem{sublemma}[theorem]{Sublemma}

\theoremstyle{definition}
\newtheorem{defn}[theorem]{Definition}
\newtheorem{example}[theorem]{Example}
\newtheorem{examples}[theorem]{Examples}
\newtheorem{remarks}[theorem]{Remarks}
\newtheorem{remark}[theorem]{Remark}
\newtheorem{algorithm}[theorem]{Algorithm}
\newtheorem{question}[theorem]{Question}
\newtheorem{problem}[theorem]{Problem}
\newtheorem{subsec}[theorem]{}
\newtheorem{clai}[theorem]{Claim}

\newtheorem{observ}[theorem]{Observation}
\newtheorem{prin}[theorem]{Principle}
\newtheorem{highlight}[theorem]{Highlight}

\def\toeq{{\stackrel{\sim}{\longrightarrow}}}

\def\into{{\hookrightarrow}}

\def\alp{{\alpha}}  \def\bet{{\beta}} \def\gam{{\gamma}}
 \def\del{{\delta}}
\def\eps{{\varepsilon}}
\def\kap{{\kappa}}                   \def\Chi{\text{X}}
\def\lam{{\lambda}}
 \def\sig{{\sigma}}  \def\vphi{{\varphi}} \def\om{{\omega}}
\def\Gam{{\Gamma}}   \def\Del{{\Delta}}
\def\Sig{{\Sigma}}   \def\Om{{\Omega}}
\def\ups{{\upsilon}}


\def\F{{\mathbb{F}}}
\def\BF{{\mathbb{F}}}
\def\BN{{\mathbb{N}}}
\def\Q{{\mathbb{Q}}}
\def\Ql{{\overline{\Q }_{\ell }}}
\def\CC{{\mathbb{C}}}
\def\R{{\mathbb R}}
\def\V{{\mathbf V}}
\def\D{{\mathbf D}}
\def\BZ{{\mathbb Z}}
\def\K{{\mathbf K}}
\def\XX{\mathbf{X}^*}
\def\xx{\mathbf{X}_*}

\def\AA{\Bbb A}
\def\BA{\mathbb A}
\def\HH{\mathbb H}
\def\PP{\Bbb P}

\def\Gm{{{\mathbb G}_{\textrm{m}}}}
\def\Gmk{{{\mathbb G}_{\textrm m,k}}}
\def\GmL{{\mathbb G_{{\textrm m},L}}}
\def\Ga{{{\mathbb G}_a}}

\def\Fb{{\overline{\F }}}
\def\Kb{{\overline K}}
\def\Yb{{\overline Y}}
\def\Xb{{\overline X}}
\def\Tb{{\overline T}}
\def\Bb{{\overline B}}
\def\Gb{{\bar{G}}}
\def\Ub{{\overline U}}
\def\Vb{{\overline V}}
\def\Hb{{\bar{H}}}
\def\kb{{\bar{k}}}

\def\Th{{\hat T}}
\def\Bh{{\hat B}}
\def\Gh{{\hat G}}

\def\cF{{\mathfrak{F}}}
\def\cC{{\mathmathcal C}}
\def\cU{{\mathmathcal U}}

\def\Xt{{\widetilde X}}
\def\Gt{{\widetilde G}}

\def\gg{{\mathfrak g}}
\def\hh{{\mathfrak h}}
\def\lie{\mathfrak a}

\def\XX{\mathfrak X}
\def\RR{\mathfrak R}
\def\NN{\mathfrak N}

\def\minus{^{-1}}

\def\GL{\textrm{GL}}            \def\Stab{\textrm{Stab}}
\def\Gal{\textrm{Gal}}          \def\Aut{\textrm{Aut\,}}
\def\Lie{\textrm{Lie\,}}        \def\Ext{\textrm{Ext}}
\def\PSL{\textrm{PSL}}          \def\SL{\textrm{SL}}
\def\loc{\textrm{loc}}
\def\coker{\textrm{coker\,}}    \def\Hom{\textrm{Hom}}
\def\im{\textrm{im\,}}           \def\int{\textrm{int}}
\def\inv{\textrm{inv}}           \def\can{\textrm{can}}
\def\id{\textrm{id}}              \def\Char{\textrm{char}}
\def\Cl{\textrm{Cl}}
\def\Sz{\textrm{Sz}}
\def\ad{\textrm{ad\,}}
\def\SU{\textrm{SU}}
\def\Sp{\textrm{Sp}}
\def\PSL{\textrm{PSL}}
\def\PSU{\textrm{PSU}}
\def\rk{\textrm{rk}}
\def\PGL{\textrm{PGL}}
\def\Ker{\textrm{Ker}}
\def\Ob{\textrm{Ob}}
\def\Var{\textrm{Var}}
\def\poSet{\textrm{poSet}}
\def\Al{\textrm{Al}}
\def\Int{\textrm{Int}}

\def\Smg{\textrm{Smg}}
\def\ISmg{\textrm{ISmg}}
\def\Ass{\textrm{Ass}}
\def\Grp{\textrm{Grp}}
\def\Com{\textrm{Com}}
\def\rank{\textrm{rank}}

\def\TAut{\textrm{TAut\,}}

\def\Bir{\textrm{Bir\,}}

\def\char{\textrm{char}}

\newcommand{\Or}{\operatorname{O}}

\def\tors{_\def{\textrm{tors}}}      \def\tor{^{\textrm{tor}}}
\def\red{^{\textrm{red}}}         \def\nt{^{\textrm{ssu}}}

\def\sss{^{\textrm{ss}}}          \def\uu{^{\textrm{u}}}

\def\mm{^{\textrm{m}}}
\def\tm{^\times}                  \def\mult{^{\textrm{mult}}}

\def\uss{^{\textrm{ssu}}}         \def\ssu{^{\textrm{ssu}}}
\def\comp{_{\textrm{c}}}
\def\ab{_{\textrm{ab}}}

\def\et{_{\textrm{\'et}}}
\def\nr{_{\textrm{nr}}}

\def\nil{_{\textrm{nil}}}
\def\sol{_{\textrm{sol}}}
\def\End{\textrm{End\,}}

\def\til{\;\widetilde{}\;}

\def\min{{}^{-1}}

\def\AGL{{\mathbb G\mathbb L}}
\def\ASL{{\mathbb S\mathbb L}}
\def\ASU{{\mathbb S\mathbb U}}
\def\AU{{\mathbb U}}


\title[Local-global invariants of finite and infinite groups]
{Local-global invariants \\ of finite and infinite groups: \\ around Burnside from another side}

\author[Kunyavski\u\i ] {Boris Kunyavski\u\i}

\address{Department of
Mathematics, Bar-Ilan University, 5290002 Ramat Gan, ISRAEL}
\email{kunyav@macs.biu.ac.il}

\begin{abstract}
This expository essay is focused on the Shafarevich--Tate set of a
group $G$. Since its introduction for a finite group by Burnside, it
has been rediscovered and redefined more than once. We discuss its
various incarnations and properties as well as relationships (some
of them conjectural) with other local-global invariants of groups.
\end{abstract}



\maketitle

\section*{Preface}
These notes focus on a local-global invariant $\Sha (G)$ of a group
$G$, its various incarnations, applications and possible
generalizations. The cyrillic $\Sha$ refers to one of the names of
this invariant, the Shafarevich--Tate set. At a glance, this text is
a bulk of definitions, vague questions and conjectures, mostly
compiled from numerous sources. However, some striking parallels
with other invariants of $G$ give a hope that something more
sensible is extractable from this eclectic material.

For a finite group $G$, the invariant we are talking about was
introduced by Burnside (naturally, under a different name) as early
as in 1911, in the second edition of his famous book \cite{Bur1}.
Soon enough, in the paper \cite{Bur2}, which was published in 1913,
he constructed the first example of $G$ with nontrivial $\Sha (G)$,
so the present article can be viewed as a modest commemoration of
the 100th anniversary of this event. Since then, $\Sha (G)$ has been
rediscovered more than once, each time revealing some new features.
The main goal of this paper is to try to put some order in its
numerous avatars and attract the attention of experts in both finite
and infinite groups, as well as in geometric group theory (and maybe
of those whose main interests lie outside group theory), to
intriguing interrelations and applications. Towards this end, an
attempt was made to make the bibliography as comprehensive as
possible, which resulted in disproportionately unbalanced structure
of the text, with a modestly sized body followed by a long tail.
However, one can hope that the reader will find here something
beyond a mere list of references.

As a precaution, it must be said that one should not be deceived by
a misleading hint in the title: the paper has nothing in common with
the well-known monograph by Kostrikin, where Burnside's name is
identified (by default, as a sort of common practice in certain
circles) with a notoriously hard group-theoretic problem (which, by
the way, also carries a distinctive local-global flavour). We are
indeed looking from another side, focusing on the development of a
far less advertised contribution of Burnside's.

\section{Main object}

\begin{defn} \label{def:main} \cite{On1}
Let $G$ be a group acting on itself by conjugation, $(g,x)\mapsto
gxg^{-1}$, and let $H^1(G,G)$ denote the first cohomology pointed set. The
set of cohomology classes becoming trivial after restriction to
every cyclic subgroup of $G$ is denoted $\Sha (G)$ and called the
Shafarevich--Tate set of $G$.
\end{defn}

\begin{defn} \label{def:rigid}
For the lack of a better term, we say that a group $G$ with one-element
Shafarevich--Tate set is $\Sha$-rigid.
\end{defn}

This term will be explained later, after clarifying relationships
with some rigidity phenomena.

\medskip

\begin{observ} \label{obs:rig}
$\Sha$-rigidity is often a crucial step to establishing important
properties of $G$, or of a whole class of groups. On the other hand,
groups with nontrivial $\Sha (G)$ often provide interesting examples
(or even allow one to refute long-standing conjectures). Some
instances will be given below.
\end{observ}

\begin{remark}
The terminology of Definition \ref{def:main} is originated in the
prototype of $\Sha (G)$, dating back to the 1950's when it appeared
in the context of a high-brow approach to diophantine equations and
has been remaining since then one of the favourite objects of
arithmetic geometers: given an algebraic group $A$ defined over a
number field $k$, $\Sha (A)$ is defined as the set of cohomology
classes $H^1(\Gamma , A(\bar k))$ (where the absolute Galois group
$\Gamma=\Gal(\bar k/k)$ acts naturally on $\bar k$-points of $A$)
that become trivial after restriction to every $\Gamma _v=\Gal (\bar
k_v/k_v)$, where $v$ runs over all places of $k$. In the purely
group-theoretic setting as above, a much more down-to-earth
description is available.

The following important remark is due to M.~Mazur (see \cite{On3}).

\end{remark}

\begin{observ} \label{rem:Mazur}
A map $f\colon G\to G$ is a cocycle if and only if the map $g\colon
G\to G$ defined by $g(x)=f(x)x$ is an endomorphism. Furthermore, $f$
is a coboundary if and only if $g$ is an inner automorphism, and the
restriction of $f$ to the cyclic subgroup generated by $x\in G$ is a
coboundary if and only if $g(x)$ is conjugate to $x$. Denote by
$\Endc (G)$ (resp. $\Autc (G)$) the set of endomorphisms (resp.
automorphisms) $g$ of $G$ such that $g(x)$ is conjugate to $x$ for
all $x\in G$.

We see that $G$ is $\Sha$-rigid if and only if it satisfies the
following condition:
\begin{equation} \label{cond:E}
\Endc (G)=\Inn (G).
\end{equation}
\end{observ}

\begin{remark} \label{rem:condE}
Endomorphisms (or automorphisms) $g$ with the property that $g(x)$
is conjugate to $x$ for all $x\in G$ appear in literature under
different names: pointwise inner, conjugating, class-preserving,
etc. In this text they will be called {\it locally inner}. Note that
any locally inner endomorphism is injective.

Condition (\ref{cond:E}) is sometimes called Property E (see, e.g.,
\cite{AKT3}).
\end{remark}

Since $\Inn (G)\subseteq \Autc (G)\subseteq \Endc (G)$, it is
convenient to subdivide the property of being $\Sha$-rigid into two
weaker ones: 1) $\Inn (G)=\Autc (G)$; 2) $\Autc (G)=\Endc (G)$. The
first property is sometimes referred to as property A (see, e.g.,
\cite{Gros}). It can be written down as $\Outc (G)=1$ meaning that
$G$ has no locally inner outer automorphisms. In this text, a group
$G$ satisfying 1) will be called A-rigid, and a group satisfying 2)
will be called B-rigid. In these terms, any $\Sha$-rigid group is
both A-rigid and B-rigid, and {\it vice versa}.

\section{Zoo of rigid groups}

\subsection{B-rigid groups} \label{sec:B-rigid}

\begin{observ} \label{prop:B-rigid}
The following groups are B-rigid:
\begin{itemize}
\item[(i)] finite;
\item[(ii)] profinite;
\item[(iii)] solvable;
\item[(iv)] cohopfian.
\end{itemize}
\end{observ}

(i) is obvious because every locally inner endomorphism is injective
and hence, as $G$ is finite, surjective. (ii) is proved in
\cite{On3}. (iii) is proved in \cite{AE}. (iv) is obvious (recall
that the property of $G$ to be cohopfian means that $G$ contains no
proper subgroups isomorphic to $G$, or, equivalently, that every
injective endomorphism of $G$ is surjective). Note that this
property is related to other rigidity properties: for example, it is
satisfied by rigid hyperbolic groups \cite{RS} (see
\cite[Theorem~4.4]{Sel} for a generalization), and irreducible
lattices in semisimple Lie groups, except for free groups (because
of their Mostow rigidity) \cite{Pr}. (Note that free groups are also
B-rigid and, moreover, $\Sha$-rigid, see Observation
\ref{prop:sha-rigid}(ii) below.) The cohopfian property also holds
for some Kleinian \cite{DP}, \cite{OP}, \cite{WZ}, 3-manifold
\cite{PV}, \cite{WaWu}, and braid groups \cite{BM}, \cite{Belli}, as
well as for some torsion-free nilpotent groups \cite{Beleg}.

\subsection{A-rigid groups} \label{sec:A-rigid}

Our first observation is obvious.

\begin{observ} \label{obs:complete}
All complete groups are A-rigid.
\end{observ}

Recall that a group $G$ is complete if it is centreless and all its
automorphisms are inner. See \cite{Rob} for a survey of finite
complete groups. As to infinite groups, typical examples of complete
groups arise as groups of automorphisms $G=\Aut(F)$, where $F$ is a
free group (or a group that is ``not so far'' from free). The cases
where $F$ is a free group, free nilpotent group of class two, or the
quotient of a free group by an appropriate characteristic subgroup
were treated, respectively, in \cite{DF1}--\cite{DF3} (for groups of
finite rank) and \cite{To1}--\cite{To3} (for groups of infinite
rank). The cases $F=\GL(n,\mathbb Z)$ ($n$ odd), $F=\PGL(2,\mathbb
Z)$ ($n\ge 2$), $F=\SL(n,\mathbb Z)$ ($n\ge 3)$ were considered in
\cite{Dy}. The cases $F=B_n$ (Artin braid group), $n\ge 4$, and
$F=\Aut(B_3)$ were established in \cite{DG}. All infinite symmetric
groups are complete \cite{DM}.

\begin{observ} \label{prop:A-rigid}
The following groups are A-rigid:

I) Finite groups:

\begin{itemize}
\item[(i)] symmetric groups \cite{OW2};
\item[(ii)] simple groups \cite{FS};
\item[(iii)] $p$-groups of order at most $p^4$ \cite{KV1};
\item[(iv)] $p$-groups having a maximal cyclic subgroup \cite{KV2};
\item[(v)]  extraspecial  \cite{KV2} and almost extraspecial $p$-groups;
\item[(vi)] $p$-groups having a cyclic subgroup of index $p^2$
\cite{KV3}, \cite{FN};
\item[(vii)] groups such that the Sylow $p$-subgroups are cyclic for odd $p$, and either
cyclic, or dihedral, or generalized quaternion for $p=2$ \cite{He1}
(see \cite{Su}, \cite{WalC} for a classification of such groups);
\item[(viii)] Blackburn groups \cite{He4}, \cite{HL};
\item[(ix)] abelian-by-cyclic groups \cite{HJ};
\item[(x)] primitive supersolvable groups \cite{La};
\item[(xi)] unitriangular matrix groups over $\mathbb F_p$ and the
quotients of their lower central series \cite{BVY};
\item[(xii)] central products of A-rigid groups \cite{KV2}.
\end{itemize}

The only new case here is that of almost extraspecial groups (see,
e.g., \cite{CT} for the definition and classification). In
particular, every such group is a central product of an extraspecial
group and a cyclic group, so the result follows from (xii).

See \cite{Ya4} for a survey and some details.

II) Infinite groups:

\begin{itemize}
\item[(i)] the absolute Galois group of $\mathbb Q$ \cite{Ik};
\item[(ii)] the absolute Galois group of  $\mathbb Q_p$ \cite{Ik}
(or, more generally, of any its finite extension (Ikeda (unpublished), \cite{JR});
\item[(iii)] non-abelian free groups \cite[Lemma 1]{Gros}; goes back to \cite{Ni};
\item[(iv)] non-abelian free profinite groups \cite{Ja};
\item[(v)] so-called pseudo-$p$-free profinite groups \cite{JR};
\item[(vi)] free nilpotent groups \cite{En};
\item[(vii)] non-abelian free solvable groups \cite{Rom};
\item[(viii)] nontrivial free products \cite{Nes1};
\item[(ix)] one-relator groups of the form $\left<a,b | [a^m,b^n]=1\right>, m, n>1$ \cite{TM};
\item[(x)] non-abelian free Burnside groups of large odd exponent \cite{Ch}, \cite{At};
\item[(xi)] fundamental groups of compact orientable surfaces \cite{Gros};
\item[(xii)] Artin braid groups $B_n$ and pure braid groups $P_n$ \cite{DG}, \cite{Nes2};
\item[(xiii)] connected compact topological groups \cite{McM};
\item[(xiv)] fundamental groups of closed surfaces with negative Euler characteristic \cite{BKZ};
\item[(xv)] non-elementary subgroups $H$ of hyperbolic
groups $G$ such that $H$ does not normalize any nontrivial finite
subgroup of $G$ \cite[Corollary~5.4]{MO};
\item[(xvi)] some groups of automorphisms and birational
automorphisms of the plane and space \cite{De1}--\cite{De4},
\cite{KS};
\item[(xvii)] unitriangular matrix groups over $\mathbb Q$ and the
quotients of their lower central series, as well unitriangular
matrix groups over $\mathbb Z$ \cite{BVY};
\item[(xviii)] all finitely generated Coxeter groups \cite{CMi}.
\end{itemize}
\end{observ}

\subsection{$\Sha$-rigid groups} \label{sec:sha-rigid}

\begin{observ} \label{prop:sha-rigid}
The following groups are $\Sha$-rigid:
\begin{itemize}
\item[(i)] groups appearing on the lists of both Observations \ref{prop:B-rigid}
and  \ref{prop:A-rigid};
\item[(ii)] free groups (\cite{OW1} for finitely generated free groups and
\cite{AKT1} in general);
\item[(iii)] groups $\SL (n,R)$, $\PSL (n,R)$ and $\GL (n,R)$ where $R$ is a
euclidean domain \cite{On4}, \cite{Wad1}, \cite{Wad2};
\item[(iv)] free products of at least two nontrivial groups \cite{AKT3};
\item[(v)] amalgamated products $A*_HB$ where $H$ is a maximal cyclic
subgroup of $A$ and $B$ \cite{AKT3};
\item[(vi)] all Fuchsian groups $G(n,r,s)$ except, possibly, triangle groups
$G(0,0,3)$ \cite{AKT3};
\item[(vii)] almost all orientable Seifert groups, except possibly $G_1(0,3)$
and $G_1(1,1)$ \cite{AKT4};
\item[(viii)] some ``polygonal'' \cite{Ki} and ``tree'' products \cite{KT};
\item[(ix)] the Cremona group of birational automorphisms of the
complex projective plane \cite{De3};
\item[(x)] all torsion-free hyperbolic groups;
\item[(xi)] right angled Artin groups \cite[Proposition~6.9 and
Remark~6.10]{Mi2}.
\end{itemize}
\end{observ}

The lists presented above look scattered, so some remarks and
questions are in order.

\begin{remark}
In cases (i)--(v), (viii)--(x), (xii), (xiv), (xv) of Observation
\ref{prop:A-rigid}.II, the groups under consideration satisfy a
stronger rigidity property: every normal automorphism (i.e., an
automorphism preserving normal subgroups) is inner (in case (iii)
this is proved in \cite{Lub}, \cite{Lue}, in case (xii) in
\cite{Nes2}, and in case (xv) in \cite{MO}; results of the latter
paper were generalized in \cite{MarMin}). Note that item (x)
provides a ``counter-example'' to the statement made in the last
paragraph of the preface. Groups in (i) are even more rigid:
according to \cite{Neu}, every automorphism of $\Gal(\mathbb Q)$ is
normal, and therefore it is inner.
\end{remark}

\begin{remark}
The connectedness assumption in (xiii) is essential in light of the
existence of finite groups that are not A-rigid.
\end{remark}

\begin{remark}
For topological groups $G$, another rigidity property is sometimes
used. In \cite{CDG} it is proved that if $G$ is a connected linear
real reductive Lie group, then every automorphism of $G$ preserving
unitary equivalence classes of unitary representations is inner. For
compact groups this property is equivalent to be locally inner so
one recovers A-rigidity, as in (xiii). Does A-rigidity hold in the
noncompact case?
\end{remark}

\begin{remark}
A-rigidity of a finitely generated group $G$, together with its
so-called conjugacy separability, implies that the group $\Out(G)$
is residually finite \cite{Gros}. This observation provides an
important source of residually finite groups, see \cite{CMi} and
references therein for details.
\end{remark}

\begin{remark}
Regarding case (xvi) of Observation \ref{prop:A-rigid}.II, let $G$
be one of the following groups: $\Aut (\mathbb A^2_{\mathbb C})$,
the affine Cremona group of the polynomial automorphisms of the
affine complex plane; $\TAut (\mathbb A^n_{\mathbb C})$, the
subgroup of the affine Cremona group consisting of the tame
polynomial automorphisms (for $n=2$ it coincides with the previous
one); $\Bir (\mathbb P^2_{\mathbb C})$, the Cremona group of
birational automorphisms of the complex projective plane. Then every
automorphism of $G$ is inner, up to composition with an automorphism
of $\mathbb C$ (see \cite{De1}, \cite{KS}, \cite{De2},
respectively). Thus the corresponding groups are all A-rigid.
Moreover, in \cite{De3} it was shown that every endomorphism of
$\Bir (\mathbb P^2_{\mathbb C})$ is a composition of an inner
automorphism with an endomorphism of $\mathbb C$, and hence $\Bir
(\mathbb P^2_{\mathbb C})$ is $\Sha$-rigid (a nontrivial
endomorphism of $\mathbb C$ cannot be locally inner). See \cite{De4}
for a survey of other rigidity properties of the Cremona group.
\end{remark}

\begin{remark}
Case (x) of Observation \ref{prop:sha-rigid} is a consequence of the
cohopfian property for freely indecomposable torsion-free hyperbolic
groups (proved in \cite{Sel}) and of \cite[Corollary~5.4]{MO}. If
the hyperbolic group is freely decomposable, then one can use
Observation \ref{prop:sha-rigid}(iv). The statement also follows
from more general results of \cite{BV}. (I thank the referee for
this remark, as well as for pointing out case (xi) of Observation
\ref{prop:sha-rigid}.)
\end{remark}

\begin{remark}
Unitriangular matrix groups over prime fields and the quotients of
their lower central series are all $\Sha$-rigid, in view of
Observations \ref{prop:B-rigid}(iii), \ref{prop:A-rigid}.I(xi) and
\ref{prop:A-rigid}.II(xvii). By \cite{BVY}, over a field that is not
prime, none of these groups is A-rigid.
\end{remark}

\begin{question}
Let $\mathcal G$ be a split simple Chevalley group over a prime
field $k$, $\mathcal B$ be a Borel subgroup of $\mathcal G$,
$\mathcal U$ be the unipotent radical of $\mathcal B$. Let
$G=\mathcal U(k)$ be the group of $k$-rational points of $\mathcal
U$. Is it $\Sha$-rigid? A similar question can be asked in the case
where $k$ is replaced with $\mathbb Z$.
\end{question}

\begin{remark} \label{Weyl}
The case of symmetric groups $\Sym (n)$ in Observation
\ref{prop:A-rigid}.I(i) is obvious because they all (except for
$n=2$ and $n=6$) fall into the class of complete groups in view of a
classical theorem of H\"older \cite{Ho}, and the exceptional cases
are easily treated separately. We put this case as a separate item
because of the following natural questions:
\end{remark}

\begin{question} \label{q:Weyl}
\begin{itemize}
\item[(i)] Are all Weyl groups $\Sha$-rigid?
\item[(ii)] What about other Coxeter groups?
\item[(iii)] What about other reflection groups?
\end{itemize}
\end{question}

A partial answer to (ii) is provided by the paper \cite{CMi},
generalizing some earlier results \cite{FH}, \cite{HRT}, see
Observation \ref{prop:A-rigid}.II(xviii). It is applicable to finite
Coxeter groups and thus answers (i) in the affirmative. (It would be
instructive to get a more elementary independent proof of (i), using
results of \cite{Ba}.) As to (iii), the case of finite complex
reflection groups, where there is a classification due to Shephard
and Todd (see, e.g., \cite{Co}) and a nice description of
automorphisms \cite{MarMic}, also looks tractable.

\begin{remark}
In view of Observations \ref{prop:A-rigid}.I(ii),
\ref{prop:A-rigid}.II(xiii) and \ref{prop:sha-rigid}(iii), it is
natural to ask what can happen if $G$ is the group of points of an
arbitrary simple linear algebraic group (or some generalization of
such). Note that the proofs of the results recorded above do not
admit straightforward generalizations (some of them are purely
computational, as in Observation \ref{prop:sha-rigid}(iii), and
others use methods of functional analysis, as in Observation
\ref{prop:A-rigid}.II(xiii), in order to deduce the needed property
from a stronger one; the proof of Observation
\ref{prop:A-rigid}.I(ii), relying on the classification of
automorphisms, can apparently be generalized).

Here are several cases where one can expect $\Sha$-rigidity and that
seem tractable.
\end{remark}

\begin{conj} \label{conj-simple}
The following groups $G$ are $\Sha$-rigid:
\begin{itemize}
\item[(i)] $G=\mathcal G(k)$, the group of $k$-points of a split simple
Chevalley group $\mathcal G$ defined over a sufficiently large field
$k$;
\item[(ii)] the same as in (i), with $k$ replaced with some ``good''
ring;
\item[(iii)] the same as in (i), with any isotropic $k$-group $\mathcal G$;
\item[(iv)] the same as in (i), with $\mathcal G$ an anisotropic group
splitting over a quadratic extension of $k$;
\item[(v)] $G\subset \mathcal G(k)$ is a ``big'' subgroup possessing
some rigidity properties in the sense of Mostow, Margulis, and
others;
\item[(vi)] $G$ is a split Kac--Moody group over a sufficiently large field $k$.
\end{itemize}
\end{conj}

Here are some comments. Let us start with (i). In this case, here is
a sketch of a possible proof of A-rigidity: according to Steinberg,
every automorphism is a composition of inner, graph, field and
diagonal automorphisms (see \cite{St}, \cite{Hum}); consider graph,
field and diagonal automorphisms separately; the first two types
should move some semisimple conjugacy class, and the third moves
some unipotent class. In case (ii), we have a Steinberg-like
classification of automorphisms \cite{Ab}, \cite{Bun}, \cite{Klj},
and can proceed as in case (i). In cases (iii) and (iv), one can use
\cite{BT} and \cite{We1}, respectively. Groups appearing in (v) were
discussed in \cite{We2}. For groups in (vi), a classification of
automorphisms is available \cite{Ca}, see also \cite{CMu1},
\cite{CMu2}; the case of finite ground fields should be similar to
\cite{FS}. Perhaps one can also treat unitary forms of Kac--Moody
groups over $\mathbb C$. As to B-rigidity, one can also use the
approach of \cite{BT} (apparently, to treat (iv)--(vi), it is to be
generalized in an appropriate way).

\begin{remark}
Trying to extend A-, B-, or $\Sha$-rigidity to other simple (or
almost simple) groups, one should not be overoptimistic. The
finitary symmetric group $\FSym (\Omega )$, where $\Omega$ is an
infinite countable set, is not B-rigid \cite{AE}. Similar arguments
can be used to show that the finitary alternating group $\FAlt
(\Omega)$ (which is simple) is not B-rigid either. The group of
automorphisms of each of these two groups is isomorphic to the
infinite symmetric group $\Sym (\Omega )$ (see, e.g., \cite{DM}), it
contains locally inner outer automorphisms (conjugation by any
element of $\Sym (\Omega )$ applied to a finitary permutation gives
the same result as conjugation by some finitary permutation), hence
none of these groups is A-rigid.
\end{remark}

\subsection{Genericity of rigid groups}

After walking around the zoo of rigid groups in the previous
sections, we would like to address the following question: are the
species we described there rare or common? How probable is to meet
any of them in wildlife? We will formulate an answer as a vague
principle.

\begin{prin} \label{prin:gen}
A reasonable property of a reasonable mathematical object lying
inside a reasonable class of objects may not hold but it will hold
at least for an object in general position (if not always), provided
the class under consideration is enlarged or restricted, if
necessary, in an appropriate way.
\end{prin}

In loose terms, this means that the probability of the event that a
randomly picked animal is a penguin, will be much higher if the
samples allowed to test are restrictively placed within the
Antarctic region.

Anyone can find lots of examples confirming this principle, looking
at his/her favourite area of mathematics. Here is an example from
number theory. Consider the following property of finite Galois
extensions $L/K$ of number fields: every element of $K$ which is a
norm everywhere locally is a norm globally (the Hasse principle).
This property fails to hold for most Galois extensions. However, if
we restrict our attention to {\it cyclic} extensions, it always
holds (Hasse). In another direction, if we allow $L/K$ to vary among
all finite extensions (not necessarily normal), the Hasse principle
holds in general position, i.e., for an extension $L/K$ of degree
$n$ such that the Galois group $\Gal (M/K)$ of the normal closure
$M$ of $L$ is the symmetric group $\Sym (n)$ \cite{VK}. More
generally, the Hasse principle for principal homogeneous spaces of
algebraic tori holds for generic maximal tori in simple algebraic
groups \cite{VK}, \cite{Kl}.

Some instances of Principle \ref{prin:gen} in the context of
$\Sha$-rigidity will be given below. Given a class of groups $G$, it
is usually a challenging conceptual task to enlarge or restrict it
in the spirit of this principle. Needless to add that to get a
meaningful mathematical statement, one has to convert expressions
such as ``general position'', ``random'', ``generic'', ``typical'',
and similar euphemisms, into a precise definition. Such a goal is
far beyond the framework of the present article, and the interested
reader is referred to relevant literature (see, e.g., \cite{Grom1},
\cite{Grom2}, \cite{Ols}, \cite{Oll}, \cite{KS1}, \cite{KS2} and
references therein) for different approaches to genericity.

Let us now go over to more concrete considerations.

\begin{examples} \label{ex:gen}
\begin{itemize}
\item[(i)]
Let $\mathcal{FP}$ denote the class of finite primitive permutation
groups. Then the class of $\Sha$-rigid groups is generic within
$\mathcal{FP}$. Indeed, according to \cite{LP}, a random element of
the symmetric group $\Sym (n)$ is contained in a primitive group
other than $\Sym(n)$ or $\Alt (n)$ with probability tending to zero
as $n\to\infty$, and it remains to refer to Observations
\ref{prop:A-rigid}.I(i) and (ii). According to the same paper
\cite{LP}, this statement remains true if one considers the class of
finite transitive permutation groups.
\item[(ii)]
Let $\mathcal{L}$ denote the class of linear groups over fields.
Then the class of $\Sha$-rigid groups is generic within
$\mathcal{L}$. This follows from  Observation
\ref{prop:sha-rigid}(ii) because a random linear group is free \cite{Ao}.
\item[(iii)]
In a similar spirit, in \cite{KS2} it was shown that ``random''
quotients of the modular group $PSL(2,\mathbb Z)$ are complete and
cohopfian, hence $\Sha$-rigid.
\item[(iv)]
In view of overwhelming genericity of torsion-free hyperbolic groups
(see, e.g., \cite{Ols}, \cite[Chapter~I]{Oll}), Observation
\ref{prop:sha-rigid}(x) shows that the class of $\Sha$-rigid groups
is generic within the class of all groups.
\end{itemize}

We finish this series of examples by a quote from \cite{KS2}:

`It seems likely that ``endomorphism rigidity'' is another general
aspect of ``randomness'': A random structure should not have any
endomorphisms except those absolutely required by the nature of the
structure.'

\end{examples}

This quote can be rephrased in the context of rigidity properties
considered above. Namely, the classes of A-, B-, and even
$\Sha$-rigid groups are much broader than most classically known
classes of rigid groups: say, there are A-rigid groups that are not
complete (and may be very far from complete, for example, free
groups) and B-rigid groups that are not cohopfian. However, Examples
\ref{ex:gen} show that $\Sha$-rigidity is another general aspect of
``randomness'': A random group should not have any endomorphisms
that behave locally as conjugations except those that behave in such
a way globally. In other words, the group-theoretic ``Hasse
principle'', introduced in Definitions \ref{def:main} and
\ref{def:rigid}, should hold generically.

\section{Properties of $\Sha (G)$}

{\it Case 1. Finite groups}

Recall that in this case $\Sha (G)$ is a group. It coincides with
$\Outc (G)$, the group of locally inner outer automorphisms.

\begin{observ} \label{sha-finite}
\begin{itemize}
\item[(i)]
The group $\Sha (G)$ is solvable \cite{Sah}.
\item[(ii)]
The group $\Sha (G)$ may be non-abelian \cite{Sah}.
\item[(iii)]
If $G$ is solvable and each chief factor is complemented, then $\Sha
(G)$ is supersolvable \cite{La}.
\item[(iv)]
If $G$ is supersolvable and its Frattini subgroup is trivial, then
$\Sha (G)$ is nilpotent \cite{La}.
\item[(v)] If $G$ is nilpotent of class $c$, then $\Sha (G)$ is
nilpotent of class at most $c-1$ \cite{Sah}.
\item[(vi)]
$\Sha (G)$ does not depend of the isoclinism class of $G$
\cite{Ya3}.
\end{itemize}
\end{observ}

\begin{remark}
Observation \ref{sha-finite}(ii) disproves an assertion made in
\cite{Bur1}. The smallest known counter-example is of order $2^{15}$
(in general, the construction gives, for each prime-power $q$, a
group of order $q^{5m}, m\ge 3$).
\end{remark}

\begin{remark} \label{p-groups}
\begin{itemize}
\item[(i)]
First examples of groups that are not $\Sha$-rigid appeared among
$p$-groups, they go back to Burnside \cite{Bur2}. The smallest group
that is not $\Sha$-rigid is of order 32 (see \cite{WalG}); a
classification of groups of order $p^5$ with nontrivial $\Sha (G)$
was obtained in \cite{Ya3} (and made more precise in \cite{Ya4}).
More examples of non-$\Sha$-rigid $p$-groups were found in
\cite{Ya1}, \cite{Ya4} among so-called Camina groups. See \cite{Ya2}
for some bounds for the order of $\Sha (G)$.
\item[(ii)]
Looking beyond $p$-groups, a classification of ``minimal'' groups
$G$ with $\Sha (G)\ne 1$ in the class of solvable groups all of
whose Sylow subgroups are abelian was obtained in \cite{He1}.
\item[(iii)]
See \cite{BVY} for some matrix counter-examples. Some general
constructions of a similar spirit can be found in \cite{Sah},
\cite{Sz}.
\end{itemize}
\end{remark}

{\it Case 2. Infinite groups}

Here much less is known.

The group $\Out_c(G)$ is finite if $G$ is hyperbolic \cite{MS1}, or,
more generally, relatively hyperbolic \cite{MS2}.

The following observation attributed to Passman (see the
introduction to \cite{Sah}) implies that $\Out_c(G)$ may be an
infinite simple group: this happens for $G=\FSym (\Omega )$, in
which case $\Out_c(G)$ is isomorphic to the quotient $\Sym
(\Omega)/\FSym (\Omega)$. (Note that this is in sharp contrast with
Observation \ref{sha-finite}(i).)

In general, from \cite{Mi1} it follows that $\Out_c(G)$ can be any
countable group. Hence the set $\Sha(G)$ can be as large as
possible. In the examples of \cite{Mi1}, $G$ is finitely generated
(in fact, 2-generated), satisfies Kazhdan's property (T) and has
exactly two conjugacy classes (including 1).

A possible conceptually interesting question is to understand the
situation within some natural classes of groups where such
pathologies do not arise. Perhaps, the first class to be considered
is that of compact groups.

\section{Shafarevich--Tate set vs. Bogomolov multiplier}

Throughout this section, unless otherwise stated, $G$ is a finite
group and $k$ is an algebraically closed field of characteritic
zero. Recall that the Bogomolov multiplier $B_0(G)$ is defined as
the subgroup of the Schur multiplier $H^2(G,\mathbb Q/\mathbb Z)$
consisting of the cocycles becoming trivial after restricting to all
abelian (or, equivalently, bicyclic) subgroup of $G$ \cite{Bo}. This
group coincides with the so-called unramified Brauer group of the
quotient $V/G$, where $V$ is a vector $k$-space equipped with a
faithful, linear, generically free action of $G$. The latter group
is an important birational invariant of $V/G$, in particular, it
equals zero whenever the variety $V/G$ is $k$-rational (or even
retract $k$-rational). It was introduced by Saltman and used in
constructing a counter-example to Noether's problem \cite{Sal}. The
Bogomolov multiplier allows one to compute this group solely in
terms of $G$.

\begin{observ} \label{ShaB}
\begin{itemize}
\item[(i)]
In many cases, the Bogomolov multiplier of $\Sha$-rigid groups is
zero. This is true at least for the groups listed in items (i)--(ix)
of Observation \ref{prop:A-rigid}.I. Moreover, for some of these
groups the corresponding varieties $V/G$ are rational (or at least,
retract rational). Indeed, in case (i) rationality is classically
known (E. Noether), it follows from the theorem on elementary
symmetric functions. This is also known in cases (iii) \cite{CK},
(iv) \cite{HuKa}, (vi) \cite{Ka1}. In case (ix) one knows retract
rationality \cite{Ka2}. The Bogomolov multiplier is zero in case
(ii) \cite{Ku}. Cases (v), (vii), (viii) were treated in \cite{KK}.
\item[(ii)]
$\Sha (G)$ is invariant under isoclinism, and so is $B_0(G)$: if
$G_1$ and $G_2$ are isoclinic  (i.e., have isomorphic quotients
$G_i/Z(G_i)$ and derived subgroups $[G_i,G_i]$, and these
isomorphisms are compatible), then we have $\Sha (G_1)\cong \Sha
(G_2)$ and $B_0(G_1)\cong B_0(G_2)$. For $\Sha (G)$ we recall
Observation \ref{sha-finite}(vi), and for $B_0(G)$ this was proved
in \cite{Mo2}; moreover, in \cite{BB} it was proved that if $G_1$
and $G_2$ are isoclinic, then the corresponding linear generically
free quotients $V/G_1$ and $V/G_2$ are stably birationally
equivalent, thus answering questions posed in \cite{HKK} in the
affirmative.
\item[(iii)]
There are examples of groups $G$ that are not $\Sha$-rigid but
$B_0(G)=0$. Such examples can be found among groups of order $p^5$:
for $p=2$ one always has $B_0(G)=0$ \cite{CHKP} but there is a group
with $\Sha (G)\ne 1$ \cite{WalG}; for every $p\ge 3$ there is an
isoclinism family (denoted $\Phi_7$) for which $\Sha (G)\ne 1$
\cite{Ya4} but $B_0(G)=0$ \cite{HKK}, \cite{Mo1}. Note, however,
that although $B_0(G)=0$ for all groups $G$ of order 32, such a
group can give rise to a homogeneous space $X=SL_n/G$ defined over a
{\it finite} field $k$ so that the unramified Brauer group of $X$ is
not zero \cite{BDH}. Thus one may pursue in giving $\Sha (G)$ some
birational flavour.
\end{itemize}
\end{observ}

One can try to proceed along this empiric line even further, asking
the following questions.

\begin{question}
\begin{itemize}
\item[(i)]
Do the groups listed in Observations \ref{prop:A-rigid}.I(x),(xi)
satisfy $B_0(G)=0$?
\item[(ii)]
Let $G=G_1*G_2$ be a central product of groups such that
$B_0(G_1)=B_0(G_2)=0$. Is it true that $B_0(G)=0$?
\item[(ii${}^{\prime}$)]
Let $G=G_1*G_2$ be a central product of groups such that the
corresponding generically free linear quotients $V_1/G_1$ and
$V_2/G_2$ are stably rational. Is it true that so is $V/G$?
\item[(iii)]
Do there exist $\Sha$-rigid groups with nonzero Bogomolov
multiplier?
\end{itemize}
\end{question}

It is a tempting task to find a conceptual explanation of the
experimental data presented above. One can try to use a hint given
in \cite{GK}, where the group $\Sha (G)$ was naturally embedded into
the so-called second lazy cohomology group $H^2_{\ell}(O_k(G))$ of
the Hopf algebra of $k$-valued functions on $G$. This object can be
viewed as a far-reaching noncommutative generalization of the Schur
multiplier of $G$ (and coincides with it in the case where $G$ is
abelian). The construction of the embedding is far from obvious:
first, the lazy cohomology is identified with the group of
equivalence classes of invariant Drinfeld twists on $k[G]$, and then
this latter group is mapped to a pointed set $\mathcal B(G)$
consisting of the pairs $(A,b)$, where $A$ runs over normal abelian
subgroups of $G$ and $b$ is a $k^*$-valued $G$-invariant
nondegenerate alternating bilinear form on the Pontryagin dual $\hat
A$. The fibre at the pointed element is then identified with $\Sha
(G)$.

Of course, this construction cannot directly be used to reveal an
eventual relationship between $\Sha (G)$ and $B_0(G)$ (say, because
the lazy cohomology is not invariant under isoclinism). However, a
more thorough exploration does not seem completely hopeless.

Moreover, one can go beyond finite groups and try to embed an
appropriate part of $\Sha (G)$ (say, the group $\Outc (G)$ of
locally inner outer automorphisms) into an appropriate cohomology
group. In the case where $G$ is compact, one can use an approach
presented in \cite{NT}.

\section{Miscellaneous applications,  ramifications, and generalizations}

As mentioned in the preface, nonrigid groups often lead to
interesting counter-examples. We list several instances of such a
phenomenon. Throughout, unless otherwise stated, $G$ is a finite
group.

\begin{observ}
\begin{itemize}
\item[(i)]
Formanek showed \cite{Fo} that if $\Sha (G)\neq 1$ and $\varphi$ is
an appropriately chosen locally inner outer automorphism of $G$,
then the semidirect product $G_1:=G\rtimes \left<\varphi\right>$
provides a counter-example to a conjecture by Roth \cite{Rot},
asserting that each irreducible representation of $G/Z(G)$ is a
subrepresentation of the conjugation representation $G/Z(G)$ on
$\mathbb C[G]$; see \cite{KM} for minimal counter-examples in the
class of $p$-groups. On the other hand, the Roth property turned out
to be related to a recently developed ``Lie theory'' and
``noncommutative differential geometry'' on finite groups
\cite{LMR}. Does there exist any direct connection between $\Sha
(G)$ and these new theories?
\item[(ii)]
The smallest non-$\Sha$-rigid group of order 32 constructed in
\cite{WalG}, as well as more complicated examples in \cite{He1},
played a crucial role in refuting some long-standing conjectures in
the theory of integral group rings, including Higman's isomorphism
problem for integral group rings of finite groups (see \cite{He2},
\cite{He4}). Some other interesting examples, based on \cite{Sah},
were obtained in \cite{PW}.
\item[(iii)]
Examples of nilpotent Lie groups with $\Sha (G)\ne 1$ were used in
\cite{GW} to construct compact Riemannian manifolds that are
isospectral but not isometric.
\item[(iv)]
Apart from $\Sha (G)$, there are other local-global invariants of a
similar flavour. One of such is related to the notion of Coleman
automorphism \cite{He3}, \cite{HeKi}: this is an automorphism of $G$
that becomes inner after restricting to each Sylow subgroup of $G$.
An argument attributed to Ph. Gille (see \cite{Pa}) shows that any
group admitting a non-inner Coleman automorphism provides an example
of a principal homogeneous space defined over a number field that
has a rational zero-cycle of degree one (=has rational points in
extensions of coprime degrees) but has no rational points. Similar
arguments were used in \cite{GGHZ} to construct non-isomorphic
curves becoming isomorphic over extensions of coprime degrees. The
existence of a principal homogeneous space of a {\it connected}
linear algebraic group over an arbitrary field with the same
property as above (with a rational zero-cycle of degree one and
without rational points) is an open problem (known as  Serre's
conjecture).
\end{itemize}
\end{observ}

\begin{remark}
As shown in the papers \cite{CMi}, \cite{BV} cited above, some
classes of $\Sha$-rigid groups $G$ are rigid in an even stronger
sense: any endomorphism of $G$ which preserves the conjugacy classes
of all elements of {\it short length} (in some suitable metric) must
be inner. This observation gives rise to a number of testability
problems of the following flavour. We say that a subset $S$ of a
$\Sha$-rigid group $G$ is {\it a test subset} if the following
holds: an endomorphism $\varphi$ of $G$ is inner if and only if
$\varphi(s)$ is conjugate to $s$ for all $s\in S$. Then one can ask
about the existence of ``small'' test subsets. Can one choose $S$
computable in some reasonable sense? finite? one-element? See
\cite{BV} for relevant discussions.
\end{remark}

Finally, one can introduce a local-global invariant more general
than $\Sha (G)$.

\begin{defn} \label{def:family}
Let $G$ be a group, let $\Gamma \leq G$ be a subgroup acting on $G$
by conjugation, and let $H^1(\Gamma ,G)$ denote the first cohomology
pointed set. Let $\mathcal F$ be a family of subgroups of $\Gamma$.
We define
$$
\Sha_{\mathcal F}(\Gamma ,G):=\ker [H^1(\Gamma ,G)\to
\prod_{\Gamma'\in\mathcal F}H^1(\Gamma' ,G)]
$$
and call it the Shafarevich--Tate set with respect to $\mathcal F$.

In particular, if $\Gamma =G$, we abbreviate $\Sha_{\mathcal
F}(G,G)$ to $\Sha_{\mathcal F}(G)$.
\end{defn}

With this terminology, if $G$ is finite then the solvable group
$\Sha (G)$ has several natural subgroups: $\Sha_{\mathcal A}(G)$,
$\Sha_{\mathcal N}(G)$, $\Sha_{\mathcal S}(G)$, where $\mathcal A$,
$\mathcal N$ and $\mathcal S$ stand for the family of all abelian,
nilpotent and Sylow subgroups of $G$ respectively. The first one is
nilpotent (of class at most 2) \cite{Da3}, the second one is abelian
\cite{Da3}, and the third one is abelian too \cite{HeKi} (earlier
Dade \cite{Da1} proved that it is solvable \cite{Da1} and then that
it is nilpotent \cite{Da4}); in all known examples the first group
is abelian too.

Using a set-up similar to Definition \ref{def:family}, D.~Segal
treated an arithmetic local-global problem of equivalence of binary
forms \cite{Seg}, and T.~Ono studied twists of hyperelliptic curves
\cite{On5}, \cite{On6}.

The reader is welcome to provide more applications and
interrelations.

\bigskip

\noindent {\it Acknowledgements}. This research was supported by the
Israel Science Foundation, grant 1207/12. The author was supported in
part by the Minerva Foundation through the Emmy Noether Research
Institute of Mathematics. This paper was mainly written during several
visits to the MPIM (Bonn) and NCTS (Taipei) in 2008--2012. Support
of these institutions is gratefully acknowledged.

The author thanks I.~Kapovich, L.~Potyagailo and N.~Vavilov for
useful discussions and the referee for thoughtful remarks.

\bigskip

\noindent {\it Bibliographical remark}. In the list of references
below, M.~Kumar in \cite{KV1}--\cite{KV3} and M.~K.~Yadav in
\cite{BVY}, \cite{Ya1}--\cite{Ya4} refer to the same person.

\end{document}